\documentclass[11pt]{article}
\usepackage{epsf,epsfig,amsfonts}
\begin{document}
\pagestyle{plain}

\def\c12{\cot(\pi u_{12})}
\def\gn#1#2{g^{(0)}_{#1#2}}
\newcommand{\intc}{\int\limits_0^{1}}
\renewcommand{\d}{\,{\rm d}}
\renewcommand{\i}{{\rm i}}
\newcommand{\del}{\partial}
\newcommand{\sgn}{\,{\rm sign}\,}
%
\newcommand{\be}{\begin{equation}}
\newcommand{\ee}{\end{equation}\noindent}
\newcommand{\bear}{\begin{eqnarray}}
\newcommand{\ear}{\end{eqnarray}\noindent}
\newcommand{\no}{\noindent}
\date{}
\renewcommand{\theequation}{\arabic{section}.\arabic{equation}}
\renewcommand{\arraystretch}{2.5}
\newcommand{\GeV}{\mbox{GeV}}
\newcommand{\cL}{\cal L}
\newcommand{\D}{\cal D}
\newcommand{\Dhalf}{{D\over 2}}
\newcommand{\Det}{{\rm Det}}
\newcommand{\Li}{{\rm Li}}
\newcommand{\PP}{\cal P}
\newcommand{\G}{{\cal G}}
\def\R{1\!\!{\rm R}}
\def\Eins{\mathord{1\hskip -1.5pt
\vrule width .5pt height 7.75pt depth -.2pt \hskip -1.2pt
\vrule width 2.5pt height .3pt depth -.05pt \hskip 1.5pt}}
\newcommand{\symb}{\mbox{symb}}
\renewcommand{\arraystretch}{2.5}
\newcommand{\slD}{\raise.15ex\hbox{$/$}\kern-.57em\hbox{$D$}}
\newcommand{\slpartial}{\raise.15ex\hbox{$/$}\kern-.57em\hbox{$\partial$}}
\newcommand{\slG}{{{\dot G}\!\!\!\! \raise.15ex\hbox {/}}}
\newcommand{\Gd}{{\dot G}}
\newcommand{\Gund}{{\underline{\dot G}}}
\newcommand{\Gdd}{{\ddot G}}
\def\GBd12{{\dot G}_{B12}}
\def\mneg{\!\!\!\!\!\!\!\!\!\!}
\def\Mneg{\!\!\!\!\!\!\!\!\!\!\!\!\!\!\!\!\!\!\!\!}
\def\non{\nonumber}
\def\beqn*{\begin{eqnarray*}}
\def\eqn*{\end{eqnarray*}}
\def\sy{\scriptscriptstyle}
\def\footstrut{\baselineskip 12pt}
\def\square{\kern1pt\vbox{\hrule height 1.2pt\hbox{\vrule width 1.2pt
   \hskip 3pt\vbox{\vskip 6pt}\hskip 3pt\vrule width 0.6pt}
   \hrule height 0.6pt}\kern1pt}
\def\slash#1{#1\!\!\!\raise.15ex\hbox {/}}
\def\dint#1{\int\!\!\!\!\!\int\limits_{\!\!#1}}
\def\bra#1{\langle #1 |}
\def\ket#1{| #1 \rangle}
\def\vev#1{\langle #1 \rangle}
\def\rightvac{\mid 0\rangle}
\def\leftvac{\langle 0\mid}
\def\dps{\displaystyle}
\def\sy{\scriptscriptstyle}
\def\half{{1\over 2}}
\def\third{{1\over3}}
\def\fourth{{1\over4}}
\def\fifth{{1\over5}}
\def\sixth{{1\over6}}
\def\seventh{{1\over7}}
\def\eigth{{1\over8}}
\def\ninth{{1\over9}}
\def\tenth{{1\over10}}
\def\pa{\partial}
\def\ddtau{{d\over d\tau}}
\def\gf{\hbox{$\gamma_5\;$}}
\def\ie{\hbox{$\textstyle{\int_1}$}}
\def\iz{\hbox{$\textstyle{\int_2}$}}
\def\id{\hbox{$\textstyle{\int_3}$}}
\def\ldop{\hbox{$\lbrace\mskip -4.5mu\mid$}}
\def\rdop{\hbox{$\mid\mskip -4.3mu\rbrace$}}
\def\eps{\epsilon}
\def\epshalf{{\epsilon\over 2}}
\def\e{\mbox{e}}
\def\g{\mbox{g}}
\def\pa{\partial}
\def\kinb{{1\over 4}\dot x^2}
\def\kinf{{1\over 2}\psi\dot\psi}
\def\expk{{\rm exp}\biggl[\,\sum_{i<j=1}^4 G_{Bij}k_i\cdot k_j\biggr]}
\def\expp{{\rm exp}\biggl[\,\sum_{i<j=1}^4 G_{Bij}p_i\cdot p_j\biggr]}
\def\expshort{{\e}^{\half G_{Bij}k_i\cdot k_j}}
\def\expabb{{\e}^{(\cdot )}}
\def\epseps#1#2{\varepsilon_{#1}\cdot \varepsilon_{#2}}
\def\epsk#1#2{\varepsilon_{#1}\cdot k_{#2}}
\def\kk#1#2{k_{#1}\cdot k_{#2}}
\def\G#1#2{G_{B#1#2}}
\def\Gp#1#2{{\dot G_{B#1#2}}}
\def\GF#1#2{G_{F#1#2}}
\def\Dab{{(x_a-x_b)}}
\def\Dsq{{({(x_a-x_b)}^2)}}
\def\lag{( -\partial^2 + V)}
\def\PITD{{(4\pi T)}^{-{D\over 2}}}
\def\4piTD{{(4\pi T)}^{-{D\over 2}}}
\def\4piT4{{(4\pi T)}^{-2}}
\def\TintmD{{\dps\int_{0}^{\infty}}{\d T\over T}\,\e^{-m^2T}
    {(4\pi T)}^{-{D\over 2}}}
\def\Tintm4{{\dps\int_{0}^{\infty}}{\d T\over T}\,\e^{-m^2T}
    {(4\pi T)}^{-2}}
\def\Tintm{{\dps\int_{0}^{\infty}}{\d T\over T}\,\e^{-m^2T}}
\def\Tint{{\dps\int_{0}^{\infty}}{\d T\over T}}
\def\pint{{\dps\int}{dp_i\over {(2\pi)}^d}}
\def\Dx{\dps\int{\cal D}x}
\def\Dy{\dps\int{\cal D}y}
\def\Dpsi{\dps\int{\cal D}\psi}
\def\Tr{{\rm Tr}\,}
\def\tr{{\rm tr}\,}
\def\sumij{\sum_{i<j}}
\def\freeexp{{\rm e}^{-\int_0^Td\tau {1\over 4}\dot x^2}}
\def\arraystretch{2.5}
\def\Ge{\mbox{GeV}}
\def\dA{\partial^2}
\def\DA{\sqsubset\!\!\!\!\sqsupset}
\def\FFdual{F\cdot\tilde F}
%
%
\def\Nset{{\mathbb{N}}}
\def\bbbr{{\mathbb{R}}}
\def\bbbone{{\mathchoice {\rm 1\mskip-4mu l} {\rm 1\mskip-4mu l}
    {\rm 1\mskip-4.5mu l} {\rm 1\mskip-5mu l}}}
\def\bbbz{{\mathbb{Z}}}
%

\pagestyle{empty}
\renewcommand{\thefootnote}{\fnsymbol{footnote}}
\vskip .4cm
\begin{center}
{\Large\bf A Quantum Field Theoretical Representation
of Euler-Zagier Sums}
\vskip1.3cm
{\large Uwe M{\"u}ller}
\\[1.5ex]
{\it Institut f{\"u}r Physik, Johannes-Gutenberg-Universit{\"a}t Mainz\\
  Staudinger-Weg 7\\
  D-55099 Mainz, Germany\\
  umueller@thep.physik.uni-mainz.de
}

\vspace{.3cm}
 {\large Christian Schubert
}
\\[1.5ex]
{\it
Laboratoire d'Annecy-le-Vieux
de Physique Th{\'e}orique LAPTH\\
Chemin de Bellevue,
BP 110\\
F-74941 Annecy-le-Vieux CEDEX,
France\\
schubert@lapp.in2p3.fr\\
} 
\vskip1.5cm

{\large \sl Published in \\
\vspace{5pt}
{\it Int. J. Math. Math. Sc.} Vol. {\bf 31}, issue 3 (2002),
127 -- 148}

\vskip1.5cm
 {\large \bf Abstract}
\end{center}

\begin{quotation}
\noindent
We establish a novel representation of 
arbitrary Euler-Zagier sums in terms of
weighted vacuum graphs.
This representation uses a toy
quantum field theory
with infinitely many propagators and 
interaction vertices.
The propagators involve Bernoulli polynomials
and Clausen functions to arbitrary orders.
The Feynman integrals of this model  
can be decomposed in terms of a vertex algebra
whose structure we investigate.
We derive a large class
of relations between multiple zeta values,
of arbitrary lengths and weights,
using only a certain set of graphical manipulations on 
Feynman diagrams.
Further uses and possible generalizations 
of the model are pointed out.

\end{quotation}
\clearpage
\renewcommand{\thefootnote}{\protect\arabic{footnote}}
\pagestyle{plain}

\setcounter{page}{1}
\setcounter{footnote}{0}

\section{Introduction}
\renewcommand{\theequation}{1.\arabic{equation}}
\setcounter{equation}{0}

The perturbative evaluation of Green's
functions in quantum field theory lead to
a class of iterated parameter integrals
whose explicit calculation becomes
very difficult beyond the first
few orders in the coupling constant expansion.
Any progress in this area of
work must be based on an intimate knowledge of
the properties of
various types of special functions such as
polylogarithms, hypergeometric functions,
and their generalizations (see e.g.\ \cite{smirnov}).

In recent years, some structure is seen to emerge
from the seemingly haphazard occurence
of those special functions as the values
of individual Feynman diagrams. Kreimer's
hypothesis \cite{kreimer}, based on
a rule of associating knots to Feynman
diagrams, allows one to predict 
from knot-theoretical considerations the  
level of transcendentality which can possibly
appear in the counterterm
coefficients of an ultraviolet divergent
diagram.
Even though it has been verified for a large number of examples
\cite{brokre96} the raison d'{\^e}tre for the
correspondence between graphs and knots remains
presently mysterious.
More recently there are indications
 that knot-theoretical concepts may
be of relevance even for the finite parts of
Feynman diagrams \cite{broadhurst98}.

The remarkably rich mathematical structures
surfacing in this correspondence make Feynman
diagrams increasingly interesting from
the pure mathematician's point of view.
The objects encountered in the calculation
of UV divergences in perturbative
quantum field theory, multiple harmonic sums,
are of considerable relevance
to number theory and other branches of
mathematics (see e.g.\ 
\cite{drinfeld,zagier,lemurakami,goncharov98,terasoma,terasoma2,furusho,racinet}).

Quantum field theory amplitudes can be
calculated in coordinate space or in momentum space.
In four-dimensional field theory the arising
integrals are normally of a similar type and
degree of difficulty. 
This is very different in the case of 
a one-dimensional quantum field theory
compactified on a circle, which will be 
considered in the present paper.
Such quantum field theories arise naturally
if one represents one-loop amplitudes
in $D$-dimensional field theory in
terms of first-quantized path integrals.
An approach to quantum field theory
along these lines
has gained some popularity in recent years
after it was discovered that
it allows one to reorganise 
ordinary field theory
amplitudes in a manner similar to
string theory amplitudes \cite{polyakov,berkos,strassler,ss3,zako}.
In this type of formalism the $D$-dimensional
space-time enters as a target space,
and amplitudes are calculated in 
terms of an auxiliary field theory
in one-dimensional parameter space.
Green's functions in parameter space
are then used 
for the evaluation of Feynman diagrams
in this one-dimensional 
`worldloop' theory.

As a simple example, let us
consider the one-loop
effective action for 
a scalar field theory with a 
${\lambda\over 3!}\phi^3$ interaction.
This effective action can
be expressed in terms of a
first-quantized path integral as follows,
\be
\Gamma[\phi]
=\half
\Tintm
\int_{x(T)=x(0)}{\cal D}x(\tau)\,
\e^{-\int_0^T\d\tau\Bigl(
\kinb
+\lambda \phi(x(\tau))
\Bigr)}.
\label{scalftpi}
\ee\no
Here $T$ is the usual Schwinger proper-time
for the particle circulating in the loop.
At fixed $T$ a path integral has to be calculated
over the space of closed loops in spacetime with
period $T$. This integral contains a zero mode
which is removed by fixing the center-of-mass
of the loop $x_0\equiv {1\over T}\int_0^T\d\tau x(\tau)$.
The reduced path integral is evaluated perturbatively
by expanding the interaction exponential and using
the parameter space Green's function
\footnote{The constant part of this Green's
function is irrelevant for the final physical results
and usually deleted from the beginning.}
\bear
G(\tau_1,\tau_2) &\equiv &
 2T
\sum_{{n=-\infty}\atop{n\ne 0}}^{\infty}
{\e^{2\pi\i n{{\tau_1-\tau_2}\over T}}
\over {(2\pi\i n)}^2}
=
\mid \tau_1-\tau_2\mid 
-{{(\tau_1-\tau_2)}^2\over T}
-{T\over 6}.
\label{defG}
\ear\no
Momentum space methods have also sometimes
been used, which in this case lead
to Fourier sum representations.
In \cite{basvan,dpsv}
a number of such sums were calculated to
provide a check on the resolution of
certain ambiguities which in the
first-quantized formalism can arise 
in curved backgrounds.
In this comparison one finds that terms given
by simple polynomial integrals in
coordinate space may, in momentum space,
correspond to non-trivial
multiple sums of the
Euler-Zagier type.

In the present work we will turn the
logic around, and use this formalism
as a tool for the systematic 
study of Euler-Zagier sums.
Euler-Zagier sums, also called multiple $\zeta$ values
or multiple harmonic series,
are defined by
\bear
\zeta(k_1,\ldots,k_m) &=&
\sum_{n_1>n_2>\cdots >n_m>0} 
{1
\over
n_1^{k_1}\cdots n_m^{k_m}
}
=
\Li_{k_1 \ldots k_m}(1,\ldots,1). \non\\
\label{defmultizeta}
\ear
They are special values of the
multidimensional polylogarithms
$\Li_{k_1,\ldots,k_m}$, defined as
\bear
\Li_{k_1\ldots k_m}(z_1,\ldots,z_m) &\equiv&
\sum_{n_1>n_2>\cdots >n_m>0} 
{z_1^{n_1}\cdots z_m^{n_m}
\over
n_1^{k_1}\cdots n_m^{k_m}
}.
\label{defmultipoly}
\ear\no
One calls $m$ the length (or depth) of such a series, and
$k_1 +k_2 + \cdots + k_m$ its level (or weight).
Sums of the type (\ref{defmultizeta}) were first
considered by Euler \cite{euler}.
Euler himself noted that
numerous relations exist between Euler-Zagier
sums. Some simple examples are the following
(all given by Euler)
\bear
\zeta(2,1) &=& \zeta (3), \label{euler1}\\
\zeta(3,1) &=& {3\over 2}\zeta (4) -\half \zeta^2 (2)
= {\pi^4\over 360}, \label{euler2}\\
\zeta(2,2) &=& \half \zeta^2 (2) -\half \zeta (4)
= {\pi^4\over 120}, \label{euler3}\\
\zeta(3,2) &=& -{11\over 2}\zeta (5) +3\zeta (2) \zeta (3),
\label{euler4}\\
\zeta (4,1) &=& 2\zeta (5) -\zeta (2)\zeta (3).
\label{euler5}
\ear\no
Further results for the length two case can be found in
\cite{tornheim,apovu}.
Systematic investigations of Euler-Zagier sums of length
higher than two have been undertaken only in recent years
\cite{raosub,hoffman92,markett,hofmoe,granville,
borgir,hoffman97,bobrbr,ohno,arakan,hofohn}.
From the point of view of physics
the study of their relations is relevant
for attempts at a classification of the possible
ultraviolet divergences in quantum field theory
\cite{brdekr,brokre96}.

We would like
to be able to represent arbitrary such sums in terms of
one-dimensional Feynman diagrams.
To achieve this goal we have to
generalize the usual worldline path integral
formalism in the following ways:

\begin{enumerate}

\item
As explained above,
the first-quantized
loop path integral is defined to run
over the space of all periodic
functions,
with the constant functions eliminated.
Here we will restrict it to one
``chiral half'' spanned by the basis
functions
$f_n(u) = \e^{2 \pi\i n u}, n = 1,2,\ldots$.
(This amounts to a complexification of spacetime.)

\item
We choose the kinetic term
of our model in such a way that arbitrary inverse powers
of derivatives will appear.

\end{enumerate}
Those purely mathematical considerations lead us to define
the ``$\zeta$-model'', a one-dimensional quantum field
theory given by the following partition function,
\bear
\lefteqn{Z(g,\lambda) = \int_{\cal H}{\cal D}x(u)\,\e^{-S},}\non\\
&&S = \int_0^1\d u_1\int_0^1\d u_2\, 
\half \bar x(u_1)\,
\Bigl(1 - \lambda 2\pi\i{\partial}^{-1}\Bigr)
x(u_2)
- \int_0^1\d u\, \e^{gx(u)+ \bar g \bar x(u)}.
\non\\
\label{zetaaction}
\ear
Here and in the following $\partial \equiv {\d\over\d u}$
denotes the ordinary derivative.
The path integral is to be performed over the Hilbert space
\be
{\cal H} =
\Big\lbrace
x(u)
\Big|
x(u)= \sum_{n=1}^{\infty}a_n\e^{2\pi\i n u},
\quad
\sum_{n=1}^{\infty} 
|a_n|^2 < \infty
\Big\rbrace.
\label{hilbert}
\ee\no
The perturbative expansion of both the kinetic and the
interaction terms for this toy model leads to the
following Feynman rules,

\begin{figure}[ht]
\begin{picture}(250,120)(10,-50)
\put (10,40){Vertices ${\mathcal V}^{p,q}$:}
\put (122,43){\circle*{3}}
\put (122,43){\vector(1,1){12}}
\put (122,43){\vector(1,-1){12}}
\put (122,43){\vector(4,1){15}}
\put (122,43){\vector(4,-1){15}}
\put (106,57){\vector(1,-1){15}}
\put (106,29){\vector(1,1){15}}
\put (122,52){\circle*{1}}
\put (122,34){\circle*{1}}
\put (126,50){\circle*{1}}
\put (126,36){\circle*{1}}
\put (118,50){\circle*{1}}
\put (118,36){\circle*{1}}
\put (180,43){$p,q = 0,1,2,\ldots$}
\put (112,20){${\bar g}^p g^q$}
\put (10,-10){Propagators:}
\put (112,-2){k}
\put (90,-9){1}
\put (100,-7){\circle*{3}}
\put (100,-7){\vector(1,0){20}}
\put (120,-7) {\line(1,0){15}}
\put (135,-7){\circle*{3}}
\put (140,-9){2}
\put (180,-10){$k =0,1,2,\ldots$}
\put (85,-30){${\lambda}^k{(2\pi\i)}^kg_{12}^{(k)}$}
\end{picture}
\caption{\label{feynmanrules}
Feynman rules of the $\zeta$-model.}
\end{figure}

\no
In section 2 we write the propagators $g_{12}^{(k)}$
explicitly in terms of Bernoulli polynomials and
Clausen functions. We show that any multiple $\zeta$ sum 
(\ref{defmultizeta}) can be represented
as a Feynman diagram in this model.
In section 3 we investigate the properties
of the elementary tree-level $n$-point integrals with
arbitrary external propagators. 
Section 4 demonstrates how one can use 
partial integrations and reality conditions
to derive a large class of relations between
multiple $\zeta$ sums. 
In section 5 we point out possible further
uses of the model, as well as some
generalizations.

\section{Basic Properties}
\renewcommand{\theequation}{2.\arabic{equation}}
\setcounter{equation}{0}

Inverting the kinetic part of our (non-local)
Lagrangian (\ref{zetaaction}), and writing $\partial^{-k}$ in
the defining basis of the Hilbert space
$\cal H$, we find
\bear
g^{(k)}_{12}
&\equiv&
g^{(k)}(u_{12})
=
\sum_{n=1}^{\infty}
{\e^{2\pi\i nu_{12}}
\over
{(2\pi\i n)}^k}.
\label{defgk}
\ear\no
If we represent the unit circle in the
complex plane, this sum will turn into
the $k$-th polylogarithm,
\bear
g^{(k)}_{12}
&=&
{1\over
{(2\pi\i)}^k}
\,\Li_k\Big({z_1\over z_2}\Big).
\label{gkpolylog} 
\ear
($u_{12} = u_1 - u_2$, $z_i = \e^{2\pi\i u_i}$).
We note the following properties of $g^{(k)}$,
\bear
\frac{\partial}{\partial u_1}g^{(k)}_{12}
=- \frac{\partial}{\partial u_2}g^{(k)}_{12}
&=& g^{(k-1)}_{12}, \label{dgk}\\
g^{(k)}_{21} &=& {(-1)}^k\bar g^{(k)}_{12}
     = g^{(k)}(1-u_{12}), \label{conjugategk}\\
g^{(k)}(0) &=& {\zeta (k) \over {(2\pi\i)}^k},
\label{gk0=zeta}\\
\int_0^1\d u_{1,2}\, g^{(k)}_{12} &=& 0.
\label{intg=0}
\ear\no
(\ref{dgk}) can be 
inverted using the explicitly known integral kernel for
inverse derivatives in this space \cite{ss3},
\bear
g_{12}^{(k+l)}&=&
\int_0^1\d u \, 
\langle u_1\mid {\partial}^{-l} \mid u\rangle\,
g^{(k)}(u-u_2),
\label{invertdelgk}\\
\langle u_1\mid {\partial}^{-n} \mid u_{2}\rangle
&=&
-{1\over n!}B_n(\vert u_{12}\vert)\,
{\rm sign}^n(u_{12})\label{delPinv1}\\
&=&
- {B_n(u_{12}) \over n!} + {u_{12}^{n-1}\over 2 (n-1)!}
({\rm sign}(u_{12})-1).\label{delPinv2}
\ear\no
Here $B_k$ denotes the $k$-th Bernoulli polynomial.
To write $g_{12}^{(k)}$ more explicitly
we split it into its real and imaginary parts,
\bear
g_{12}^{(k)} &=&
\bar g_{12}^{(k)} + \i\, \hat g_{12}^{(k)}.
\label{splitgk}
\ear\no
Using the standard integral representation
of the polylogarithm
\bear
\Li_k(z) &=&
{(-1)^{k-1}\over (k-1)!}\,z
\int_0^1 \d x\, {\ln^{k-1}(x)\over 1-xz},
\label{intrepLik}
\ear
it is then easy to show that
\pagebreak
\bear
g_{12}^{(0)}
&=&
\half
\Bigl(
\delta (u_{12})-1+\i\,\c12
\Bigr),
\label{g0}\\
g_{12}^{(1)}
&=&
\half
\Bigl(
\half {\rm sign} (u_{12}) - u_{12} 
+{\i\over \pi}
\ln |2\sin (\pi u_{12})|
\Bigr),
\label{g1}\\
g_{12}^{(2)}
&=&
\fourth
\Bigl(
|u_{12}|
-u_{12}^2
-{1\over 6}
+
{\i\over \pi^2}
\int_0^1
\d\xi\,
{\ln \xi \sin(2\pi u_{12})
\over
{1-2\xi\cos(2\pi u_{12})
+\xi^2}}
\Bigr),
\non\\
\label{g2}\\
g_{12}^{(3)}
&=&
-{1\over 24}({\rm sign}(u_{12})-2u_{12})(|u_{12}|-u_{12}^2)
\non\\&&
+{\i\over 16\pi^3}
\int_0^1{\d\xi\over\xi}\ln\xi\ln(1-2\xi\cos (2\pi u_{12}) +\xi^2),
\non\\
\label{g3}\\
g_{12}^{(k\ne 0\,\, {\rm even})}
&=&
-\half {1\over k!}B_k(|u_{12}|)
\non\\&&
+{\i\over {(2\pi)}^k}{{(-1)}^{{k\over 2}+1}\over (k-1)!}
\int_0^1\d\xi\, \ln^{k-1}\xi
{\sin (2\pi u_{12})\over 1-2\xi\cos (2\pi u_{12})+\xi^2},\non\\
\label{gkeven}\\
g_{12}^{(k\ne 1\,\,{\rm odd})}
&=&
-\half {1\over k!}B_k(|u_{12}|){\rm sign}(u_{12})\non\\
&&
+{\i\over 2{(2\pi)}^k}{{(-1)}^{{k+1\over 2}}\over (k-2)!}
\int_0^1\!{\d\xi\over\xi}\, {\ln}^{k-2}\xi
\ln(1-2\xi\cos (2\pi u_{12})+\xi^2),\non\\
\label{gkodd}
\ear
($u_{12}\in [-1,1]$).
Note also that 
the imaginary part of $g_{12}^{(k)}$ is,
up to a normalization 
factor, identical with the $k$-th Clausen function
${\rm Cl}_k(2\pi u_{12})$ (see, e.g., \cite{lewin}).
Only the real parts are present in the calculation of worldline
path integrals in standard field theory. In particular,
note that $\bar g^{(2)}$ is, up to a conventional factor
of $4$, identical with the function $G$ introduced in
(\ref{defG}) (for $T=1$).
Only the
imaginary parts are capable of producing $\zeta (n)$'s
with $n$ odd.

\begin{figure}[ht]
\begin{center}
\begin{picture}(0,0)%
\epsfig{file=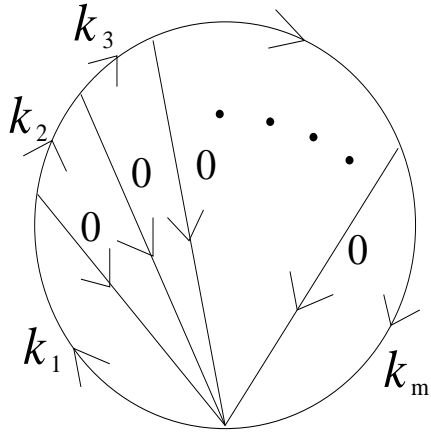,height=160pt,width=160pt}%
\end{picture}%
\setlength{\unitlength}{0.00087500in}%
\begingroup\makeatletter\ifx\SetFigFont\undefined
\def\x#1#2#3#4#5#6#7\relax{\def\x{#1#2#3#4#5#6}}%
\expandafter\x\fmtname xxxxxx\relax \def\y{splain}%
\ifx\x\y   
\gdef\SetFigFont#1#2#3{%
  \ifnum #1<17\tiny\else \ifnum #1<20\small\else
  \ifnum #1<24\normalsize\else \ifnum #1<29\large\else
  \ifnum #1<34\Large\else \ifnum #1<41\LARGE\else
     \huge\fi\fi\fi\fi\fi\fi
  \csname #3\endcsname}%
\else
\gdef\SetFigFont#1#2#3{\begingroup
  \count@#1\relax \ifnum 25<\count@\count@25\fi
  \def\x{\endgroup\@setsize\SetFigFont{#2pt}}%
  \expandafter\x
    \csname \romannumeral\the\count@ pt\expandafter\endcsname
    \csname @\romannumeral\the\count@ pt\endcsname
  \csname #3\endcsname}%
\fi
\fi\endgroup
\begin{picture}(3834,2736)(120,-1050)
\end{picture}
\caption{\label{sea shell}
`Sea shell' diagram representing the general Euler-Zagier sum.}
\end{center}
\end{figure}

A vacuum diagram in our model involves a multiple
integral on the unit circle with an integrand which is 
a product
of propagators (\ref{defgk}). (One of the integrations
is redundant due to translation invariance.)
The result of the integrations can obviously be
decomposed as a sum of multiple $\zeta$ values.
Our model thus defines a map from the set of
weighted vacuum graphs to the (vector space
generated by) the Euler-Zagier sums.
This map is, moreover, easily seen to be surjective;
the reader will immediately convince
herself that with the above Feynman 
rules the following ``sea shell''
diagram (fig.\ \ref{sea shell}) 
evaluates to
\bear
\lambda^{\sum k_i}(g\bar g)^{2m-1}
\,
\zeta(k_1,\ldots,k_m).
\label{valuesea shell}
\ear\no
The $u$-integrations produce $\delta$ functions
for ``momentum conservation'' at every vertex, while the
Fourier sums of the inserted $g^{(0)}$ propagators
yield Heaviside step functions leading to the desired
ordering of the remaining sums.

\no

\section{The elementary vertex integrals}
\renewcommand{\theequation}{3.\arabic{equation}}
\setcounter{equation}{0}

Let us now start on an investigation of the
properties of the Feynman integrals in
$x$($= u$)-space. An obvious first step is to
consider the folding of the elementary vertices
with arbitrary sets of propagators. 
We denote the elementary vertex integral by
\bear
I_{k_1\ldots k_p}^{l_1\ldots l_q}(u_1,\ldots,u_{p+q})
&\equiv&
\int_0^1\d u \, g^{(k_1)}(u_1-u)\cdots g^{(k_p)}(u_p-u)
\non\\&&\times            
g^{(l_1)}(u-u_{p+1})\cdots g^{(l_q)}(u-u_{p+q}).
\ear\no
We note that it has the following obvious properties,
\bear
I_{k_1\ldots k_p}^{l_1\ldots l_q} &=& 0 \qquad 
\mbox{if $p=0$  or $q=0$},
\label{allinallout}\\
\bar I_{k_1\ldots k_p}^{l_1\ldots l_q} &=&
{(-1)}^{\sum_{i=1}^p k_i + \sum_{j=1}^q l_j}
I^{k_1\ldots k_p}_{l_1\ldots l_q}.
\label{Icc}
\ear\no

\subsection{Two-vertex integral}
\label{sec:2vertex}

By construction two-point vertices can be
integrated out trivially,
\bear
\int_0^1\d u_3 \, g_{13}^{(k)}g_{32}^{(l)}
=
g_{12}^{(k+l)}.
\label{int2vert}
\ear

\subsection{Three-vertex integral}
\label{sec:3vertex}

The evaluation of vertex integrals is complicated by singularities
which can appear due to the presence of the cotangent function
in $g^{(0)}$.
Integrals involving $\cot(\pi u_{12})$ need to be performed
using the principal value prescription.
One way of calculating them is to transform them into
complex contour integrals via the substitution $z=\exp(2\pi\i u)$,
\begin{equation}
  \intc\d u\prod_{k=1}^n\cot(\pi(u-u_k))
  =\i^n\oint\frac{\d z}{2\pi\i z}\prod_{k=1}^n\frac{z+z_k}{z-z_k}.
  \label{contour-ints}
\end{equation}
Those 
are evaluated by means of residues (the poles on the contour give
half values due to the principal value prescription)
\footnote{We remark that, if one expresses the propagators
of the sea shell diagram via (\ref{gkpolylog}) and
(\ref{intrepLik}), and then calculates the $u_i$-integrals
by means of residues, then one arrives precisely at
Kontsevich's integral representation \cite{kassel} for the
multiple $\zeta$ sum.}. The first few
integrals are
\begin{eqnarray}
  \intc\d u\,\cot(\pi(u-u_1))&=&0,
  \label{imag-part-int1}\\
  \intc\d u\prod_{k=1}^2\cot(\pi(u-u_k))&=&
  \delta(u_1-u_2)-1,
  \label{imag-part-int2}\\
  \intc\d u\prod_{k=1}^3\cot(\pi(u-u_k))&=&
  \delta(u_1-u_2)\cot(\pi(u_2-u_3))+\nonumber\\[-3mm]
  &+&\delta(u_1-u_3)\cot(\pi(u_1-u_2))+\nonumber\\
  &+&\delta(u_2-u_3)\cot(\pi(u_2-u_1)).
  \label{imag-part-int3}
\end{eqnarray}\no
Alternatively one can also calculate those integrals
recursively using, under the integral, the following
identity,
\begin{eqnarray}
  \lefteqn{
  \cot(\pi u_{12})\cot(\pi u_{13})+
  \cot(\pi u_{21})\cot(\pi u_{23})+}\nonumber\\
  &&+\cot(\pi u_{31})\cot(\pi u_{32})=
  -1+\delta(u_{12})\delta(u_{13}).
  \label{cot-identity}
\end{eqnarray}\no
This identity will be of further use later on.
With these results the vertex integral of three propagators $g^{(0)}$
is 
\begin{eqnarray}
  I_{00}^0 (u_1,u_2,u_3) 
&=&\intc\d u\, g^{(0)}(u_1-u)g^{(0)}(u_2-u)g^{(0)}(u-u_3)
  \nonumber\\
  &=&
  \frac{1}{4}\left[
  \delta(u_{13})\delta(u_{23})-\delta(u_{13})-\delta(u_{23})+1
  \right.
  \nonumber\\
  &&+(\delta(u_{13})-1)\i\cot(\pi u_{23} )\nonumber\\
  &&+(\delta(u_{23})-1)\i\cot(\pi u_{13})\nonumber\\
  &&\left.
  -\cot(\pi u_{23})\cot(\pi u_{13})\right],
  \label{3-vertex}
\end{eqnarray}
and can be identified as
\begin{equation}
 I_{00}^0 (u_1,u_2,u_3)  
=g^{(0)}_{13}g^{(0)}_{23}.
  \label{3-vertex:a}
\end{equation}
Applying the identity (\ref{cot-identity}) gives also
\begin{eqnarray}
  I_{00}^0 (u_1,u_2,u_3) 
&=&
  g^{(0)}_{12}g^{(0)}_{23}-
  g^{(0)}_{12}g^{(0)}_{13}- \frac{1}{2}g^{(0)}_{13}\nonumber\\
  &=&
  g^{(0)}_{21}g^{(0)}_{13}-
  g^{(0)}_{21}g^{(0)}_{23}-
  \frac{1}{2}g^{(0)}_{23}\nonumber\\
  &=&
  g^{(0)}_{12}g^{(0)}_{23}+
  g^{(0)}_{21}g^{(0)}_{13}-
  \frac{1}{2}\delta_{12}g^{(0)}_{13}.
  \label{3-vertex:b}
\end{eqnarray}
Folding of eq.\ (\ref{3-vertex:a}) with
$\bra{u_{1'}}\del^{-n}\ket{u_1}$ and $\bra{u_{2'}}\del^{-m}\ket{u_2}$
leads to
\begin{eqnarray}
I_{nm}^0 (u_1,u_2,u_3) &=&
  \intc\d u\, g^{(n)}(u_1-u)g^{(m)}(u_2-u)g^{(0)}(u-u_3)
  \nonumber\\
  &=&g^{(n)}_{13}g^{(m)}_{23},
  \label{3-vertex:12-complete}
\end{eqnarray}
whereas the forms (\ref{3-vertex:b}) can be folded with
$\bra{u_3}\del^{-k}\ket{u_{3'}}$
\begin{eqnarray}
I_{00}^k (u_1,u_2,u_3) &=&
  \intc\d u\, g^{(0)}(u_1-u)g^{(0)}(u_2-u)g^{(k)}(u-u_3)
  \nonumber\\
  &=&
  g^{(0)}_{12}g^{(k)}_{23}-
  g^{(0)}_{12}g^{(k)}_{13}
  -\frac{1}{2}g^{(k)}_{13}\nonumber\\
  &=&
  g^{(0)}_{21}g^{(k)}_{13}-
  g^{(0)}_{21}g^{(k)}_{23}-
  \frac{1}{2}g^{(k)}_{23}\nonumber\\
  &=&
  g^{(0)}_{12}g^{(k)}_{23}+
  g^{(0)}_{21}g^{(k)}_{13}-
  \frac{1}{2}\delta_{12}g^{(k)}_{13}. 
  \label{3-vertex:3-complete}
\end{eqnarray}\no
Nontrivial is the case
\bear
I_{0k}^l (u_1,u_2,u_3) &=&
\int_0^1 \d u\, g^{(0)}(u - u_1)g^{(k)}(u-u_2)g^{(l)}(u_3-u)
\non\\
&=&
\Li_{lk}(z_{31},z_{12}).
\label{3-vertex:nontrivial}
\ear\no
The general case of all three indices different from zero
can be reduced to this case by partial integrations, e.g.
\bear
I_{11}^1 (u_1,u_2,u_3)
&=&
\int_0^1 \d u\, g^{(1)}(u_{41})g^{(1)}(u_{42})g^{(1)}(u_{34})
\non\\
&=&
\int_0^1 \d u\, g^{(1)}(u_{41})g^{(0)}(u_{42})g^{(2)}(u_{34})
\non\\&&
+\int_0^1 \d u\, g^{(0)}(u_{41})g^{(1)}(u_{42})g^{(2)}(u_{34})
\non\\
&=&
\Li_{21}(z_{31},z_{12})
+
\Li_{21}(z_{32},z_{21}).
\label{examppartint}
\ear\no

\subsection{A three -- point relation}
\label{sec:3rel}

\begin{figure}[ht]
\begin{center}
\begin{picture}(0,0)%
\epsfig{file=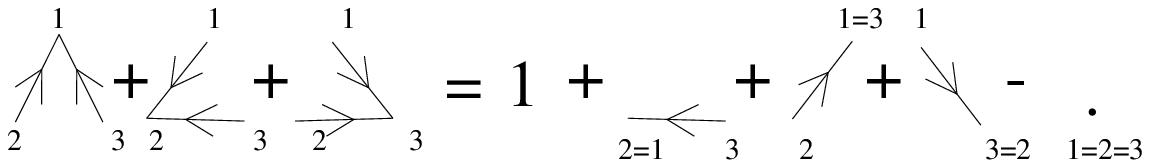}%
\end{picture}%
\setlength{\unitlength}{0.00087500in}%
\begingroup\makeatletter\ifx\SetFigFont\undefined
\def\x#1#2#3#4#5#6#7\relax{\def\x{#1#2#3#4#5#6}}%
\expandafter\x\fmtname xxxxxx\relax \def\y{splain}%
\ifx\x\y   
\gdef\SetFigFont#1#2#3{%
  \ifnum #1<17\tiny\else \ifnum #1<20\small\else
  \ifnum #1<24\normalsize\else \ifnum #1<29\large\else
  \ifnum #1<34\Large\else \ifnum #1<41\LARGE\else
     \huge\fi\fi\fi\fi\fi\fi
  \csname #3\endcsname}%
\else
\gdef\SetFigFont#1#2#3{\begingroup
  \count@#1\relax \ifnum 25<\count@\count@25\fi
  \def\x{\endgroup\@setsize\SetFigFont{#2pt}}%
  \expandafter\x
    \csname \romannumeral\the\count@ pt\expandafter\endcsname
    \csname @\romannumeral\the\count@ pt\endcsname
  \csname #3\endcsname}%
\fi
\fi\endgroup
\begin{picture}(5474,936)(120,-550)
\end{picture}
\caption{\label{3pointrel}
Three-point relation.}
\end{center}
\end{figure}

Beyond the three-point case the systematic investigation of
the elementary vertex integrals becomes cumbersome.
We will be satisfied here to note that, 
from the representation (\ref{g0}) for
$g^{(0)}$ and the identity (\ref{cot-identity}) 
one easily derives the following pair of (complex conjugate)
three-point identities,
\bear
\gn21\gn31 + \gn12\gn32 + \gn13\gn23
&=&
1 + \delta_{12} \gn32 + \delta_{31}\gn21
+\delta_{23}\gn13 
-\delta_{12}\delta_{13},
\non\\
\gn12\gn13 + \gn21\gn23 + \gn31\gn32
&=&
1 + \delta_{12} \gn23 + \delta_{31}\gn12
+\delta_{23}\gn31 
-\delta_{12}\delta_{13}.
\non\\
\label{3pointid}
\ear
We may represent the first
identity graphically as in fig. \ref{3pointrel}.

\no
The second identity has all arrows reversed.
Those identities can be used to transform any vertex integral involving
two $g^{(0)}$'s which are either both ingoing or both outgoing.

By iteration of this three-point identity one can contruct
an analogous identity for an arbitrary number of points.
The formulas are rather cumbersome 
and will not be given here.

\section{Derivation of Multiple $\zeta$ Relations}
\renewcommand{\theequation}{4.\arabic{equation}}
\setcounter{equation}{0}

We will now show how can use the formalism developed above for
deriving a large class of 
multiple $\zeta$ relations by simple manipulations
on graphs, with no need to ever explicitly write down
sums. In those manipulations we will make use of
the following elements:

\begin{enumerate}

\item
The triviality of the real part of $g^{(0)}$,
$\bar g^{(0)}_{12} = \half (\delta_{12} -1)$
(see (\ref{g0})).

\item
Partial integrations.

\item
The three-point relations eq.\ (\ref{3pointid}).

\item
The vanishing of diagrams containing a vertex with 
only ingoing or only outgoing propagators.

\item
The two-vertex integration formula eq.\ (\ref{int2vert}).

\end{enumerate}

At intermediate steps sometimes ill-defined sums
will appear such as $\zeta (1)$. Those will always
cancel out in the final results. The appearance
of divergent sums could be easily avoided using
a regularization such as in \cite{borgir} but we
will not do so here.

\subsection{Length two}
\label{length2}

\begin{figure}[ht]
\begin{center}
\begin{picture}(0,0)%
\epsfig{file=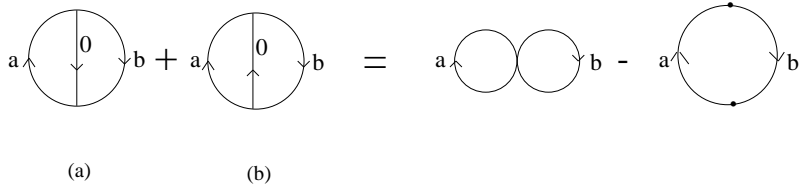}%
\end{picture}%
\setlength{\unitlength}{0.00087500in}%
\begingroup\makeatletter\ifx\SetFigFont\undefined
\def\x#1#2#3#4#5#6#7\relax{\def\x{#1#2#3#4#5#6}}%
\expandafter\x\fmtname xxxxxx\relax \def\y{splain}%
\ifx\x\y   
\gdef\SetFigFont#1#2#3{%
  \ifnum #1<17\tiny\else \ifnum #1<20\small\else
  \ifnum #1<24\normalsize\else \ifnum #1<29\large\else
  \ifnum #1<34\Large\else \ifnum #1<41\LARGE\else
     \huge\fi\fi\fi\fi\fi\fi
  \csname #3\endcsname}%
\else
\gdef\SetFigFont#1#2#3{\begingroup
  \count@#1\relax \ifnum 25<\count@\count@25\fi
  \def\x{\endgroup\@setsize\SetFigFont{#2pt}}%
  \expandafter\x
    \csname \romannumeral\the\count@ pt\expandafter\endcsname
    \csname @\romannumeral\the\count@ pt\endcsname
  \csname #3\endcsname}%
\fi
\fi\endgroup
\begin{picture}(5264,1636)(120,1250)
\end{picture}
\vspace{30pt}
\caption{\label{reflect}
Diagrammatic representation of the reflection identity.}
\end{center}
\end{figure}

Let us begin at length two. The simplest way of arriving
at identities is provided by the first one of the points
listed above. Consider the sea shell diagram at length
two (fig.\ \ref{reflect}a), representing
$\zeta (a,b)$, and the same diagram with the middle
propagator reversed
(fig.\ \ref{reflect}b),
representing $\zeta (b,a)$.
(We will generally disregard the coupling constant factors
in the following.)  

\no
Adding up both diagrams we
can replace the middle propagator by twice its
real part 
(fig.\ \ref{reflect}). Since
$\bar g^{(0)}_{12} = \half (\delta_{12} -1)$
this diagram can then be replaced by the sum of the two diagrams 
shown in the right hand side of fig.\ \ref{reflect}.
Using eq.(\ref{int2vert}) on the rightmost one of those
we obtain the identity
\bear
\zeta (a,b) + \zeta(b,a) &=& \zeta(a)\zeta(b) - \zeta (a+b).
\label{idreflect}
\ear\no
This identity is well-known \cite{hoffman92,borgir}, and has been named
``reflection formula'' in \cite{borgir}.

\begin{figure}[ht]
  \begin{center}
    \epsfig{file=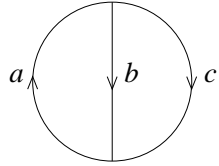}
    \caption{Auxiliary diagrams for multiple $\zeta$ relations of
      length 2.}
    \label{fig:halfmoon}
  \end{center}
\end{figure}

\begin{figure}[ht]
  \begin{center}
    \epsfig{file=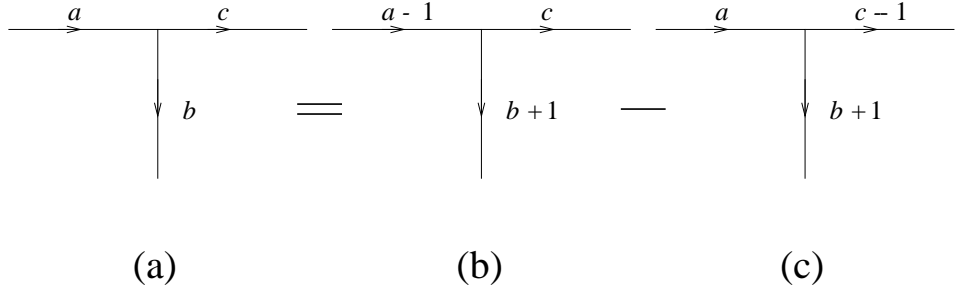,width=\textwidth}
    \caption{Integration by parts.}
    \label{fig:part-int}
  \end{center}
\end{figure}

Another way of obtaining identities is partial integration. Instead of
considering the sea shell diagram (fig.\ \ref{reflect}a), which
represents $\zeta(a,b)$, let us consider the more general diagram (fig.\ 
\ref{fig:halfmoon}) that represents a number $G_{a,b,c}=G_{a,c,b}$.
For $a=0$, $b=0$ or $c=0$, it can be represented by multiple $\zeta$
functions,
\begin{equation}
  \label{G3-zeta}
  G_{0,b,c}=\zeta(b)\zeta(c),\quad
  G_{a,b,0}=G_{a,0,b}=\zeta(a,b).
\end{equation}
To see the first of these identities one adds the same diagram with
the left propagator reversed, and uses elements one and four
from the above list.
For the later discussion we refer to these diagrams as zeta diagrams in
contrast to the other, non-zeta diagrams. If both $a$ and $c$ are
greater than $1$, we use integration by parts at the upper vertex
according to fig.\ \ref{fig:part-int} and obtain
\begin{equation}
  \label{part-int2}
  G_{a,b,c}=G_{a-1,b+1,c}-G_{a,b+1,c-1}.
\end{equation}
This can be repeated until either $a=0$ or $c=0$
\begin{eqnarray}
  \label{G-solution}
  G_{a,b,c}&=&\sum_{n=1}^c(-1)^{c+n}{a+c-n-1\choose a-1}G_{0,a+b+c-n,n}+
  \nonumber\\
  &&{}+\sum_{n=1}^a(-1)^{c}{a+c-n-1\choose c-1}G_{n,a+b+c-n,0}.
\end{eqnarray}
The right side of this equation can be translated immediately into
$\zeta$ values by (\ref{G3-zeta}). 
For $b=0$,
we obtain the relation
\begin{eqnarray}
  \label{zeta-len2}
  \zeta(a,b)=(-1)^b&\Big[&
    \sum_{n=1}^b(-1)^n{a+b-n-1\choose a-1}\zeta(n)\zeta(a+b-n)
    +\nonumber\\
    &&{}+\sum_{n=1}^a {a+b-n-1\choose b-1}\zeta(n,a+b-n)
\Bigr].
\end{eqnarray}
The divergent $\zeta(1)$ and
$\zeta(1,a)$ appear here always in the combination 
\bear
\zeta(1)\zeta(a+b-1)-\zeta(1,a+b-1),
\label{combzetadiv}
\ear\no
which can be
reexpressed in terms of convergent sums 
by the reflection identity (\ref{idreflect}). 
In this way we arrive at
\begin{eqnarray}
  \label{zeta-len2final}
  \zeta(a,b)=(-1)^b&\Big[&
    \sum_{n=2}^b(-1)^n{a+b-n-1\choose a-1}\zeta(n)\zeta(a+b-n)
    +\nonumber\\
    &&{}+\sum_{n=2}^a
{a+b-n-1\choose b-1}\zeta(n,a+b-n)-
    \nonumber\\
    &&{}-{a+b-2\choose a-1}\left(\zeta(a+b)+\zeta(a+b-1,1)\right)
  \Big].
\end{eqnarray}
Here on the left side we must assume $a>1$ if regularisation
is to be avoided.

\subsection{Length three}
\label{length3}

\begin{figure}[ht]
\begin{center}
\begin{picture}(0,0)%
\epsfig{file=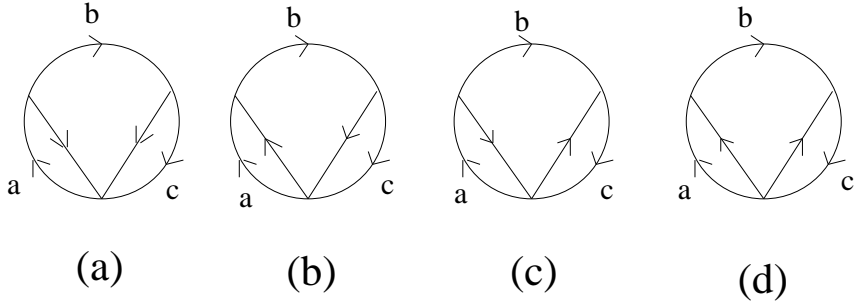}%
\end{picture}%
\setlength{\unitlength}{0.00087500in}%
\begingroup\makeatletter\ifx\SetFigFont\undefined
\def\x#1#2#3#4#5#6#7\relax{\def\x{#1#2#3#4#5#6}}%
\expandafter\x\fmtname xxxxxx\relax \def\y{splain}%
\ifx\x\y   
\gdef\SetFigFont#1#2#3{%
  \ifnum #1<17\tiny\else \ifnum #1<20\small\else
  \ifnum #1<24\normalsize\else \ifnum #1<29\large\else
  \ifnum #1<34\Large\else \ifnum #1<41\LARGE\else
     \huge\fi\fi\fi\fi\fi\fi
  \csname #3\endcsname}%
\else
\gdef\SetFigFont#1#2#3{\begingroup
  \count@#1\relax \ifnum 25<\count@\count@25\fi
  \def\x{\endgroup\@setsize\SetFigFont{#2pt}}%
  \expandafter\x
    \csname \romannumeral\the\count@ pt\expandafter\endcsname
    \csname @\romannumeral\the\count@ pt\endcsname
  \csname #3\endcsname}%
\fi
\fi\endgroup
\begin{picture}(5054,2036)(160,1050)
\end{picture}
\vspace{20pt}
\caption{\label{3shellg0}
Diagrams related to $\zeta(a,b,c)$. The propagators without
label are $g^{(0)}$ propagators.}
\end{center}
\end{figure}

\no
Proceeding to sums of length three, let us again begin
by exploiting the triviality of $\bar g^{(0)}$. 
Since there are two
$g^{(0)}$ propagators  
we have now several possibilities.
Fig.\ \ref{3shellg0} shows the ``standard'' sea shell
diagram, representing $\zeta (a,b,c)$,
as well as the diagrams related to it by a change of
direction of one or both of the $g^{(0)}$ propagators.

\vspace{-11pt}
\noindent
Those diagrams represent the quantities

\bear
(a) &:& \zeta (a,b,c), \non\\
(b) &:& \zeta(b,a,c) + \zeta(b,c,a) + \zeta (b,a+c), \non\\
(c) &:& \zeta(a,c,b) + \zeta(c,a,b) + \zeta (a+c,b), \non\\
(d) &:& \zeta(c,b,a). \non\\
\label{idreflect3}
\ear
\vspace{-2pt}

\no
Adding those diagrams in pairs to create $\bar g^{(0)}$'s
one obtains the following four identities,
\bear
\zeta(a,b,c)+\zeta(b,a,c)+\zeta(b,c,a)
+ \zeta(b,a+c) +\zeta(a+b,c) -\zeta(a)\zeta(b,c)\!\!&=&\!\! 0,
\non\\
\label{idreflect31}\\
\zeta(a,b,c)+\zeta(a,c,b)+\zeta(c,a,b)+\zeta(a+c,b)
+\zeta(a,b+c)-\zeta(c)\zeta(a,b)\!\! &=&\!\! 0,
\non\\
\label{idreflect32}\\
\zeta(b,a,c)+\zeta(b,c,a)+\zeta(c,b,a) 
+\zeta(b,a+c)
+\zeta(b+c,a)-\zeta(c)\zeta(b,a) \!\!&=&\!\! 0,
\non\\
\label{idreflect33}\\
\zeta(a,c,b)+\zeta(c,a,b)+\zeta(c,b,a)+\zeta(a+c,b)
+\zeta(c,a+b)-\zeta(a)\zeta(c,b) \!\!&=&\!\! 0.
\non\\
\label{idreflect34}
\ear\no
These identites generalize the reflection identity
(\ref{zeta-len2}). Let us call them
``permutation identities''.

\begin{figure}[ht]
\begin{center}
\begin{picture}(0,0)%
\epsfig{file=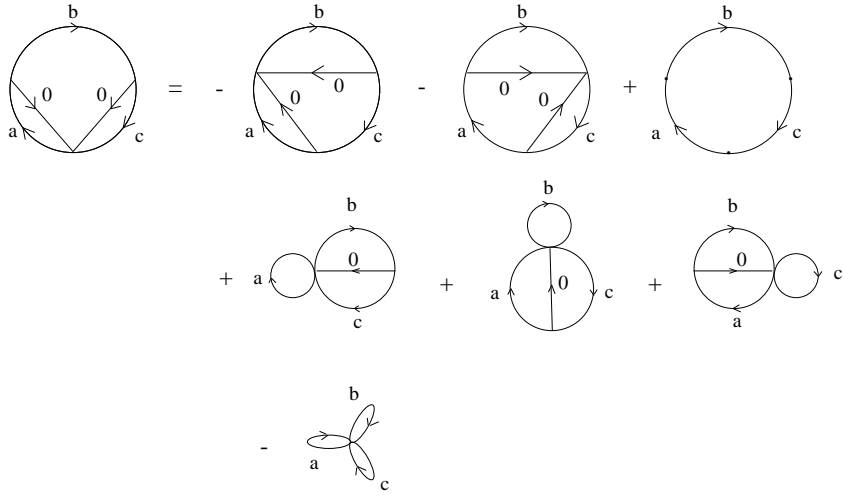}%
\end{picture}%
\setlength{\unitlength}{0.00087500in}%
\begingroup\makeatletter\ifx\SetFigFont\undefined
\def\x#1#2#3#4#5#6#7\relax{\def\x{#1#2#3#4#5#6}}%
\expandafter\x\fmtname xxxxxx\relax \def\y{splain}%
\ifx\x\y   
\gdef\SetFigFont#1#2#3{%
  \ifnum #1<17\tiny\else \ifnum #1<20\small\else
  \ifnum #1<24\normalsize\else \ifnum #1<29\large\else
  \ifnum #1<34\Large\else \ifnum #1<41\LARGE\else
     \huge\fi\fi\fi\fi\fi\fi
  \csname #3\endcsname}%
\else
\gdef\SetFigFont#1#2#3{\begingroup
  \count@#1\relax \ifnum 25<\count@\count@25\fi
  \def\x{\endgroup\@setsize\SetFigFont{#2pt}}%
  \expandafter\x
    \csname \romannumeral\the\count@ pt\expandafter\endcsname
    \csname @\romannumeral\the\count@ pt\endcsname
  \csname #3\endcsname}%
\fi
\fi\endgroup
\begin{picture}(5054,3336)(160,2050)
\end{picture}
\vspace{20pt}
\caption{\label{3shell3rel}
Diagrammatic identity derived from the three-point identity.}
\end{center}
\end{figure}

At length three we can also already make use of the
three-point identities (\ref{3pointid}).
We consider again diagram \ref{3shellg0}a, and
apply the first one of the three-point relations to
the two $g^{(0)}$ propagators running into the
``root'' vertex. The result is the
diagrammatic identity 
shown in fig.\ \ref{3shell3rel}.

\no
It can be translated term by term into the following
$\zeta$ identity,
\bear
\zeta(a,b,c) &=&
-\zeta(b,c,a)-\zeta(c,a,b)
+\zeta(a+b+c)
+\zeta(a) \zeta(b,c) \non\\&&
+ \zeta(b)\zeta(c,a)
+\zeta(c)\zeta(a,b) 
 -\zeta(a)\zeta(b)\zeta(c).
\label{id3shell3rel}
\ear\no
The identities (\ref{idreflect32}) through (\ref{idreflect34}) can be
obtained by applying the permutations $b\to a\to c\to b$,
$a\leftrightarrow c$, $b\leftrightarrow c$ on
(\ref{idreflect31}). Taking the latter identity and applying all
permutations, we obtain 6 identities which can be written as
\begin{equation}
  \label{id-matrix}
  M\vec z=\vec a,
\end{equation}
where $\vec z=(\zeta(a,b,c),\zeta(a,c,b),\ldots,\zeta(c,b,a))$ (all
permutations of the arguments) and $\vec a$ is a vector which contains
only $\zeta$ values of length 1 and 2. The rank of the coefficient
matrix $M$ is 4, i.e.\ four zeta values in $\vec z$ can be expressed
by the other two (and lower-length zeta values), e.g.\ by
$\zeta(a,b,c)$ and $\zeta(a,c,b)$.  Taking into account also relation
(\ref{id3shell3rel}) the coefficient matrix $M$ gets more rows but
the rank does not change. Therefore this relation is up to
lower-length identities not independent from (\ref{idreflect31}).
When two of the arguments $a$, $b$, $c$ coincide, say $b=c$, $\vec z$
becomes $(\zeta(a,b,b),\zeta(b,a,b),\zeta(b,b,a))$ and the rank of $M$
reduces to 2. For all arguments coinciding the rank of $M$ is 1.

There is another possibility to derive an identity from diagram
(\ref{3shellg0}d) which generalizes immediately to $\zeta$ values
of larger length. Every $g^{(0)}$ propagator in (\ref{3shellg0}d) can
be replaced by the reverted propagator according to
\begin{equation}
  \label{g0revers}
  g^{(0)}_{12}=-g^{(0)}_{21}+\delta_{12}-1.
\end{equation}
This leads to a sum of diagrams where all ocurring $g^{(0)}$
propagators are directed towards the ``root'' vertex. The resulting
identity
\begin{eqnarray}
  \label{refl-id}
  \zeta(c,b,a)&=&\zeta(a,b,c)-\zeta(a,b)\zeta(c)+\zeta(a,b+c)+
  \nonumber\\
  &&{}+\zeta(a)\left(-\zeta(b,c)+\zeta(b)\zeta(c)-\zeta(b+c)\right)+
  \nonumber\\
  &&{}+\zeta(a+b,c)-\zeta(a+b)\zeta(c)+\zeta(a+b+c)
\end{eqnarray}
can be obtained also by subtracting (\ref{idreflect31}) from
(\ref{idreflect34}) and applying appropriately the length-2 identity
(\ref{idreflect}), thus it is not a new identity.

\begin{figure}[ht]
  \begin{center}
    \epsfig{file=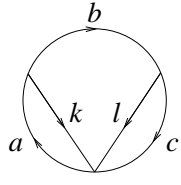}
    \caption{Auxiliary diagrams for partial integrations at length 3. }
    \label{fig:3shellgen}
  \end{center}
\end{figure}

\begin{figure}[ht]
  \begin{center}
    \epsfig{file=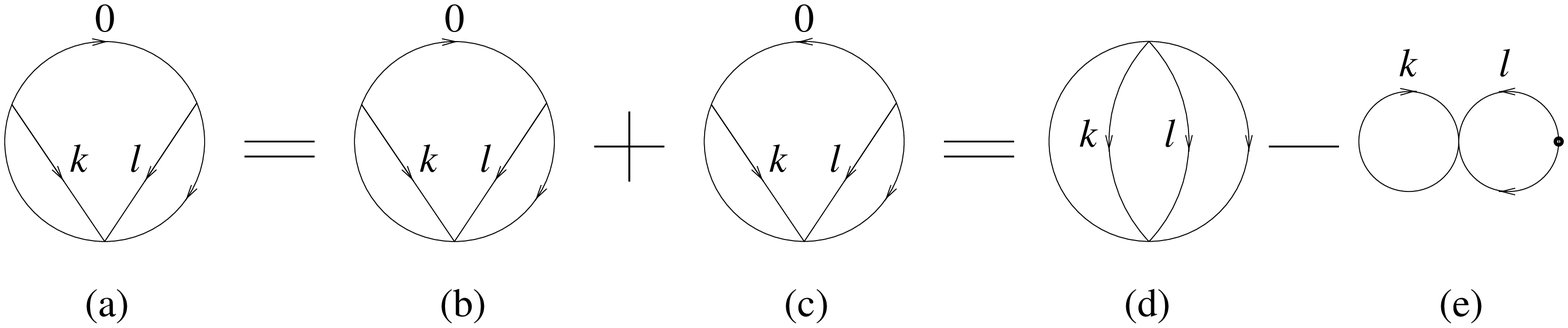,width=\textwidth}
    \caption{Identity which proves that $k$ and $l$ can be exchanged 
      when the upper propagator is zero. Diagram (c) is zero because
      no propagator goes into the upper right vertex. Diagram (e) is
      zero for the same reason at the right vertex.}
    \label{fig:part-int-exchange}
  \end{center}
\end{figure}

Partial integrations yield additional identities. As for $\zeta$
values of length 2, we consider the more general diagrams fig.\
\ref{fig:3shellgen} which evaluate to the numbers
$G_{a,k,b,l,c}=G_{a,k,b,c,l}$.
Some of these numbers can be identified with $\zeta$ values:
\begin{equation}
  \label{G5-zeta}
  G_{a,0,b,0,c}=G_{a,0,b,c,0}=\zeta(a,b,c),\quad
  G_{0,k,b,0,c}=G_{0,k,b,c,0}=\zeta(k)\zeta(b,c).
\end{equation}
We start with $\zeta(a,b,c)=G_{a,0,b,0,c}$. First, we integrate by parts
at the upper right vertex until $b=0$ or $c=0$. The combinatorics is
the same as in eq.\ (\ref{G-solution}). The terms with $c=0$ can be
identified with $\zeta$ values. In the terms where $b=0$, the two
inner propagators can be exchanged according to the identity depicted
in fig.\ \ref{fig:part-int-exchange}.
After this exchange we can integrate by parts at the upper left vertex
until terms with $a=0$ or $k=0$ arise. In terms of $\zeta$ values,
we obtain
\begin{eqnarray}
  \label{zeta-len3}
  \lefteqn{\zeta(a,b,c)=}\nonumber\\
  &=&(-1)^b\sum_{n=1}^c\Biggl[\sum_{m=1}^a
  {a+b+c-m-n-1\choose b-1,a-m,c-n}
  \zeta(m,a+b+c-m-n,n)+\nonumber\\  
  &&{}\hspace{1.5cm}+\sum_{m=1}^{b+c-n}(-1)^m{b+c-n-1\choose b-1}
  {a+b+c-m-n-1\choose a-1}\times\nonumber\\
  &&\hspace{5.5cm}\times\zeta(m)\zeta(a+b+c-m-n,n)\Biggr]+\nonumber\\
  &&{}+\sum_{n=1}^b(-1)^c{b+c-n-1\choose c-1}\zeta(a,n,b+c-n),
\end{eqnarray}
where on the left hand side we have again to require that $a>1$.
The divergent terms for $m=1$ appear always in a combination which can
be eliminated by identity (\ref{idreflect32}),
\begin{equation}
  \label{zeta3div-eliminate}
  \zeta(1)\zeta(a,b)-\zeta(1,a,b)=
  \zeta(a,b,1)+\zeta(a,1,b)+\zeta(a+1,b)+\zeta(a,b+1).
\end{equation}
Alternatively one may after the first step, in the terms with $b=0$,
rather than interchanging the two inner propagators, 
those labelled 
``k'' and ``l'' in fig.\ 9, instead interchange propagators 
``k'' and ``c''. Proceeding in the same way as before one
obtains the following identity,
\bear
\zeta (a,b,c) &=&
{(-1)}^c \sum_{n=1}^b
{b+c-n-1\choose c-1}
\zeta(a,n,b+c-n)
\non\\
&&
+ {(-1)}^c \sum_{n=1}^c\sum_{m=1}^a
{b+c-n-1\choose b-1}
{a-m+n-1\choose n-1}
\non\\
&&
\times\zeta(m,a-m+n,b-n+c)
\non\\
&&
+ \sum_{n=1}^c\sum_{m=1}^n
{(-1)}^{c-m}
{a-m+n-1\choose a-1}
{b-n+c-1\choose b-1}
\non\\
&&
\times\zeta(m)\zeta(a-m+n,b-n+c).
\label{zeta3altern}
\ear\no
Here again the terms involving a $\zeta(1,\ldots )$ appear
only in the combination (\ref{idreflect32}) and thus
can be removed without the need for regularisation.

\subsection{Arbitrary length}
\label{length4+}

All three different procedures which we have used for constructing
multiple $\zeta$ identities -- partial integrations, reversion of
$g^{(0)}$ propagators, and the use of the three-point identity --
can be generalized to the arbitrary length case without difficulty.

\begin{figure}[htb]
\epsfig{file=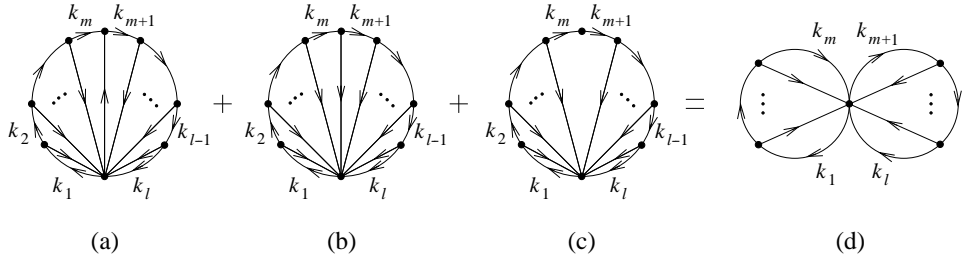,width=\textwidth}
\caption{Diagrammatic representation of permutation identities for 
arbitrary length. The propagators without label are $g^{(0)}$
  propagators.}
\label{refl-id-gen}
\end{figure}
The generalisation of the permutation identities to an arbitrary length
is based on finding representations of sea shell diagrams with
reverted inner $g^{(0)}$ propagators, which are added in order to use
the two--point relation (\ref{g0revers}). Let us consider first the
simplest case with one reverted propagator (fig.\ \ref{refl-id-gen}). 
Diagrams \ref{refl-id-gen}.b--d can be expressed
immediately by $\zeta$ values. Diagrams \ref{refl-id-gen}.a--c represent
the series
\begin{equation}
  \label{refl-id-gen-abc}
  \sum_{\raisebox{0mm}[1ex]{
      \begin{array}[t]{c}\scriptstyle n_1,\ldots,n_l\\[-8mm]
        \scriptstyle n_1>\ldots>n_m>0,\\[-8mm]
        \scriptstyle n_{m+1}>\ldots>n_l>0\end{array}}}^\infty
 \frac{1}{n_1^{k_1}n_2^{k_2}\ldots n_l^{k_l}}.
\end{equation}
By decomposing the summation range into regions where the
$n_\lambda$'s are completely ordered, we obtain a sum of $\zeta$
values which can be constructed as follows.

Let us denote by
$M_{a_1,a_2,a_{12}}\subset\{\{1\},\{2\},\{1,2\}\}^{a_1+a_2+a_{12}}$
the set of all $(a_1+a_2+a_{12})$-tuples consisting of $a_1$ elements
$\{1\}$, $a_2$ elements $\{2\}$ and $a_{12}$ elements $\{1,2\}$.
Further, we define a map ($l\equiv m+m'$)
\begin{eqnarray}
  \label{refl-id-gen-map}
  \rho_{m,m'}^a\colon
  \Nset^{l}\times M_{m-a,m'-a,a}&\to&\Nset^{l-a},
  \nonumber\\
  (k_1,\ldots,k_m;k_{m+1},\ldots,k_{l};m_1,\ldots,m_{l-a})&\mapsto&
  \nonumber\\
  &&\hspace{-2.5cm}(b_1+b_1',\ldots,b_{l-a}+b_{l-a}'),
\end{eqnarray}
where $(b_1,\ldots,b_{l-a})$ results from $(m_1,\ldots,m_{l-a})$ by
replacing $\{2\}$ with $0$, and $\{1\}$ and $\{1,2\}$ with $k_1$,
\ldots, $k_m$ (in this order); $(b_1',\ldots,b_{l-a}')$ results from
$(m_1,\ldots,m_{l-a})$ by replacing $\{1\}$ with $0$, and $\{2\}$ and
$\{1,2\}$ with $k_{m+1}$, \ldots, $k_l$ (in this order). For example,
\begin{equation}
  \label{rho-example}
  \rho_{2,2}^1(k_1,k_2;k_3,k_4;\{1\},\{1,2\},\{2\})=(k_1,k_2+k_3,k_4).
\end{equation}\no
With these definitions, the identity depicted in fig.\
\ref{refl-id-gen} can be written as
\begin{eqnarray}
  \label{refl-id-gen-result}
  \sum_{a=0}^{\min(m,l-m)}\sum_{\xi\in M_{m-a,l-m-a,a}}\hspace{-5mm}
  \zeta(\rho_{m,l-m}^a(k_1,\ldots,k_m;k_{m+1},\ldots,k_l;\xi))&=&\\
  &&\hspace{-3.5cm}\zeta(k_1,\ldots,k_m)\,\zeta(k_{m+1},\ldots,k_l).
  \nonumber
\end{eqnarray}

The generalization to the cases with more than one reverted propagator
should be obvious. But we note, that in these cases the $\zeta$ values of
maximal length appear in combinations which can be constructed also
from the identities (\ref{refl-id-gen-result}). Thus we conjecture
that all identities which are based on the reordering of summation
ranges are generated by (\ref{refl-id-gen-result}).

Imitating the considerations on page \pageref{id-matrix}, we can write
all permutation identities (\ref{refl-id-gen-result}) for fixed $l$
(but varying $m$), where the arguments of the length-$l$ $\zeta$
values are taken from the set $\{k_1,\ldots,k_l\}$ in all possible
orderings, in the form (\ref{id-matrix}), where now the components of
$\vec z$ consist of all different $\zeta$ values which result from
$\zeta(k_1,\ldots,k_l)$ by permutations of the arguments. $\vec a$
contains only $\zeta$ values of length less than $l$. Assuming $k_1$,
\ldots, $k_l$ mutually different, we found that the coefficient matrix
$M$ has rank 18 for $l=4$ and rank 96 for $l=5$. These results suggest
that for length $l$ the rank is $l!-(l-1)!$ and that the permutation
identities suffice to express the $l!$ considered $\zeta$ values of
length $l$ by the subset of $(l-1)!$ $\zeta$ values where one of the
arguments is held fixed at a certain position.

To generalize the partial integration procedure from length three
to length $m$, we can proceed in various ways. For example, we
can simply iterate the above
exchange of inner propagators. In the first step one 
applies the same partial integration as in the length three
case to the rightmost vertex of the sea shell diagram
(fig.\ \ref{sea shell}). For those terms in the result where
$k_{m-1}=0$ one makes use of this new $g^{(0)}$ propagator
to interchange the adjacent inner propagators. Then
one performs partial integrations on the next-to-rightmost
vertex until either $k_{m-2}=0$ or the next inner propagator
becomes the zero propagator. In the first case the procedure
continues with another interchange of inner propagators.
After maximally $m-1$ such propagator interchanges,
proceeding from right to left, the final
step is reached, which is again the same as in the length
three case.
 
In the first step we have the same ambiguity as before.
Depending on its resolution we arrive at a generalization
of either (\ref{zeta3div-eliminate}) or
(\ref{zeta3altern}).
We give here the formula generalizing eq.(\ref{zeta3altern}),

\bear
\zeta(k_1,\ldots,k_m)
&=&
(-1)^{k_m}
\sum_{n_{m-1}=1}^{k_{m-1}}
{k_{m-1}-n_{m-1}+n_m-1\choose k_m-1}
\non\\
&&\quad\times
\zeta(k_1,\ldots,k_{m-2},n_{m-1},k_{m-1}-n_{m-1}+n_m)
\non\\
&&\hspace{-80pt}
+
(-1)^{k_m}
\sum_{n_{m-1}=1}^{n_{m}}
\sum_{n_{m-2}=1}^{k_{m-2}}
{k_{m-1}-n_{m-1}+n_m-1\choose k_{m-1}-1}
{k_{m-2}-n_{m-2}+n_{m-1}-1\choose n_{m-1}-1}
\non\\
&&\hspace{-50pt}\times
\zeta(k_1,\ldots,k_{m-3},n_{m-2},k_{m-2}-n_{m-2}+n_{m-1},
k_{m-1}-n_{m-1}+n_m)
\non\\
&&\hspace{-80pt}
+
(-1)^{k_m}
\sum_{n_{m-1}=1}^{n_{m}}
\sum_{n_{m-2}=1}^{n_{m-1}}
\sum_{n_{m-3}=1}^{k_{m-3}}
{k_{m-1}-n_{m-1}+n_m-1\choose k_{m-1}-1}
{k_{m-2}-n_{m-2}+n_{m-1}-1\choose k_{m-2}-1}
\non\\
&&\hspace{-50pt}\times
{k_{m-3}-n_{m-3}+n_{m-2}-1\choose n_{m-2}-1}
\zeta(k_1,\ldots,k_{m-4},n_{m-3},k_{m-3}-n_{m-3}+n_{m-2},
\ldots)
\non\\
&&\hspace{-20pt}
\vdots
\non\\
&&\hspace{-80pt}
+ (-1)^{k_{m}}
\sum_{n_{m-1}=1}^{n_{m}}
\sum_{n_{m-2}=1}^{n_{m-1}}
\ldots
\sum_{n_2=1}^{n_3}
\sum_{n_1=1}^{k_1}
\prod_{n=2}^{m-1}
{k_n-n_n+n_{n+1}-1\choose k_n-1}
{k_1-n_1+n_2-1\choose n_2 -1}
\non\\
&&\hspace{-50pt}\times
\zeta(n_1,k_1-n_1+n_2,\ldots,k_{m-1}-n_{m-1}+n_m)
\non\\
&&\hspace{-80pt}
+
\sum_{n_{m-1}=1}^{n_{m}}
\sum_{n_{m-2}=1}^{n_{m-1}}
\ldots
\sum_{n_1=1}^{n_2}
(-1)^{{k_m}-n_1}
\prod_{n=1}^{m-1}
{k_n-n_n+n_{n+1}-1\choose k_n-1}
\non\\
&&\hspace{-50pt}\times
\zeta(n_1)\zeta(k_1-n_1+n_2,k_2-n_2+n_3,\ldots,k_{m-1}-n_{m-1}+n_m),
\non\\
\label{idbig}
\ear\no
where $n_m \equiv k_m$. 
Let us also give the special case $k_m=1$ of this formula
which is particularly simple,
\begin{eqnarray}
  \label{len-n-id}
  \zeta(k_1,\ldots,k_{m-1},1)&=&\zeta(1)\,\zeta(k_1,\ldots,k_{m-1})-\\
  &&\hspace{-7.5ex}{}-\sum_{\kappa=1}^{m-1}
\sum_{n_\kappa=1}^{k_\kappa}
  \zeta(k_1,\ldots,k_{\kappa-1},k_\kappa+1-n_\kappa,n_\kappa,
  k_{\kappa+1},\ldots,k_{m-1}). \nonumber
\end{eqnarray}
The terms involving $\zeta(1,\ldots)$ can be removed by means of a
special case of eq.\ (\ref{refl-id-gen-result}), namely
\begin{eqnarray}
  \label{refl-id-special}
  \lefteqn{\zeta(1)\,\zeta(k_1,\ldots,k_{m-1})
    -\zeta(1,k_1,\ldots,k_{m-1})=}\\
  &&\hspace{3cm}\sum_{\kappa=1}^{m-1}\big[
    \zeta(k_1,\ldots,k_{\kappa-1},k_\kappa+1,k_{\kappa+1},\ldots,k_{m-1})+
    \nonumber\\
    &&\hspace{3cm}{}+
    \zeta(k_1,\ldots,k_\kappa,1,k_{\kappa+1},\ldots,k_{m-1})\big].
  \nonumber
\end{eqnarray}

\begin{figure}[htbp]
  \begin{center}
    \epsfig{file=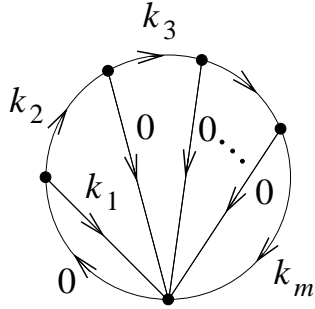,height=4cm}
    \caption{Modified sea shell diagram.}
    \label{fig:modified-seashell}
  \end{center}
\end{figure}

Up to here, the partial integrations were sequentially applied
starting at the rightmost vertex and ending at the leftmost one. An
interesting alternative is to do the opposite. Consider fig.\ 
\ref{fig:modified-seashell} (page \pageref{fig:modified-seashell}). As
in the derivation of the other partial-integration identities, this
diagram represents $\zeta(k_1)\zeta(k_2,\ldots,k_m)$. On the other
hand, we can use partial integrations at the leftmost vertex until we
have only terms where either $k_1$ or $k_2$ became $0$. If $k_1=0$
then the resulting diagram represents a multiple $\zeta$ value. If
$k_2=0$ then we exchange the adjacent inner propagators and repeat the
whole procedure at the next-to-leftmost vertex. This continues until
the rightmost vertex is reached,
where after the partial integrations all terms
represent multiple $\zeta$ values. The result is
\begin{eqnarray}
  \label{idbignice}
  \lefteqn{\zeta(k_1)\zeta(k_2,\ldots,k_m)=}\nonumber\\
  &&\sum_{n_1=1}^{k_2}{k_2+k_1-n_1-1\choose k_1-1}
\,\zeta(k_2+k_1-n_1,n_1,k_3,\ldots,k_m)+ \nonumber\\
  &&\sum_{\kappa=2}^{m-1}
  \sum_{n_1=1}^{n_0}\sum_{n_2=1}^{n_1}\ldots
  \sum_{n_{\kappa-1}=1}^{n_{\kappa-2}}\sum_{n_\kappa=1}^{k_{\kappa+1}}
  \prod_{\lambda=2}^\kappa
  {k_\lambda+n_{\lambda-2}-n_{\lambda-1}-1 \choose k_\lambda-1}
  \times\nonumber\\
  &&{k_{\kappa+1}+n_{\kappa-1}-n_{\kappa}-1 \choose n_{\kappa-1}-1}
  \times\nonumber\\
  &&\zeta(k_2+n_0-n_1,k_3+n_1-n_2,\ldots,k_{\kappa+1}+n_{\kappa-1}-n_\kappa,
  n_\kappa,k_{\kappa+2},\ldots,k_m)+\nonumber\\
  &&\sum_{n_1=1}^{n_0}\sum_{n_2=1}^{n_1}\ldots
  \sum_{n_{m-1}=1}^{n_{m-2}}
  \prod_{\lambda=2}^m
  {k_\lambda+n_{\lambda-2}-n_{\lambda-1}-1 \choose k_\lambda-1}
  \times\nonumber\\
  &&\zeta(k_2+n_0-n_1,k_3+n_1-n_2,\ldots,k_{m}+n_{m-2}-n_{m-1},
  n_{m-1}),  
\end{eqnarray}
($n_0\equiv k_1$). The right hand side contains no
divergent terms when $k_1,k_2 \ge 2$. For $k_1=1, k_2\ge 2$
it becomes finite when combined with 
(\ref{refl-id-special}).

From our derivation clearly one would expect eqs.\ (\ref{idbig})
and (\ref{idbignice}) to be equivalent. And indeed, it is a matter of 
pure combinatorics to show that (\ref{idbig}) becomes trivially
fulfilled
if (\ref{idbignice}) is used on the right hand side. Nevertheless,
we suspect that the form (\ref{idbignice}) may be more useful
for the application of these formulas to the problem of
constructing a minimal basis of independent multiple $\zeta$ sums.

\begin{figure}[htbp]
  \begin{center}
    \epsfig{file=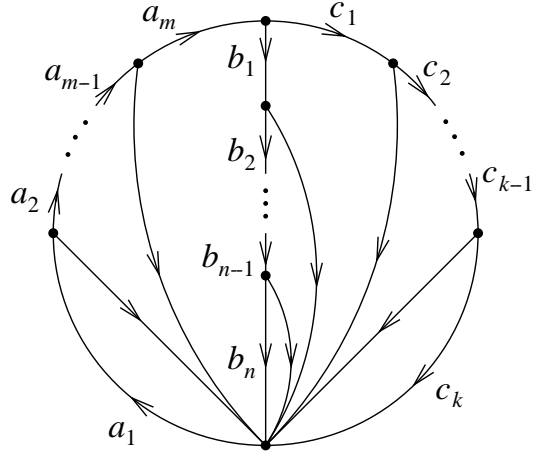,height=6cm}
    \caption{``Peacock'' diagram. The propagators without label are $g^{(0)}$
      propagators.}
    \label{fig:shuffle}
  \end{center}
\end{figure}
Eq.\ (\ref{idbignice}) is a special case of a class of identities
derived from (weight-length) shuffle algebras. In order to derive that
whole class, consider Fig.\ \ref{fig:shuffle}. It evaluates to a
number which we denote by
\begin{equation}
  \label{eq:shuffle-graph}
  Z(a_1,\ldots,a_m|b_1,\ldots,b_n|c_1,\ldots,c_k)=
  Z(a_1,\ldots,a_m|c_1,\ldots,c_k|b_1,\ldots,b_n).
\end{equation}
Partial integrations at the top vertex yield, similarly to
(\ref{G-solution}),
\begin{eqnarray}
  \label{eq:shuffle-part-int}
  \lefteqn{Z(\ldots,a_m,0|b_1,b_2,\ldots|c_1,c_2,\ldots)=}\nonumber\\
  &&\sum_{\nu=1}^{b_1}{b_1+c_1-\nu-1\choose c_1-1}
  Z(\ldots,a_m,b_1+c_1-\nu|\nu,b_2,\ldots|0,c_2,\ldots)\nonumber\\
  &&{}+\sum_{\nu=1}^{c_1}{b_1+c_1-\nu-1\choose b_1-1}
  Z(\ldots,a_m,b_1+c_1-\nu|0,b_2,\ldots|\nu,c_2,\ldots).
\nonumber\\
\end{eqnarray}
Considerations like in Fig.\ \ref{fig:part-int-exchange} show
\begin{eqnarray}
  \label{eq:shuffle-exchange}
  Z(a_1,\ldots,a_m|0,b_1,\ldots|c_1,\ldots)&=&
  Z(a_1,\ldots,a_m|c_1,\ldots|0,b_1,\ldots)\nonumber\\
  &=&Z(a_1,\ldots,a_m,0|b_1,\ldots|c_1,\ldots).\nonumber\\
\end{eqnarray}
Starting with
\begin{equation}
  \label{eq:shuffle-zeta1}
  \zeta(a_1,\ldots,a_m)\,\zeta(b_1,\ldots,b_n)=
  Z(0|a_1,\ldots,a_m|b_1,\ldots,b_n)
\end{equation}
and applying continuedly eqs. (\ref{eq:shuffle-part-int}) and
(\ref{eq:shuffle-exchange}), we end up with terms of the form
\begin{eqnarray}
  \label{eq:shuffle-zeta2}
  Z(a_1,\ldots,a_m|0|b_1,\ldots,b_n)&=&
  Z(a_1,\ldots,a_m|b_1,\ldots,b_n|0)\nonumber\\
  &=&\zeta(a_1,\ldots,a_m,b_1,\ldots,b_n).
\end{eqnarray}
The resulting identities have the same form as the shuffle identities
in the literature, but here advantageously the binomial coefficients
in (\ref{eq:shuffle-part-int}) explicitly encode (in part) the
combinatorics implicit in the shuffle algebra.

Additional multiple $\zeta$ identities can be derived
using the three-point identity (\ref{3pointid}). 
At low lengths and levels it turns out that the multiple
$\zeta$ identities obtained in this way are not independent
from the set of equations generated by the propagator
reversions and partial integrations. Whether this property
holds true in general we do not know.

\section{Discussion}
\renewcommand{\theequation}{5.\arabic{equation}}
\setcounter{equation}{0}

In the present work we have established a novel
representation of multiple $\zeta$ sums
in terms of Feynman diagrams in a
$1+0$ dimensional quantum field theory.
We
demonstrated the usefulness of this
representation
for the
derivation of identities
between such sums.  
The encoding 
into Feynman diagrams proposed here provides
a very convenient book-keeping device for
certain formal manipulations performed on such sums,
as our examples should have amply demonstrated.

Concerning the novelty of the identities derived here,
the length-two identities presented in section 4.1 are,
of course, well-known.
At length three, the identities derived by the partial
integration procedure, (\ref{zeta3div-eliminate}),
(\ref{zeta3altern}), are similar, and presumably
equivalent, to the `decomposition equations' derived in \cite{borgir}
by explicit series manipulations (their eq.(1)).
Similarly, their `permutation equations' (2)
coincide with our eq.(\ref{idreflect32}).
Eq.(\ref{id3shell3rel}) is contained as
a special case in Theorem 2.2 of \cite{hoffman92}.
However, we have not been 
able to locate in the literature an exact equivalent of our
length $m$ identities (\ref{idbig}),(\ref{idbignice})
\footnote{The special case obtained by combining
eqs.(\ref{len-n-id}) and (\ref{refl-id-special}) is Theorem
5.1 of \cite{hoffman92}.}. 
The only identities available
for arbitrary lengths and levels are those based on the
`shuffle algebra'
\cite{bobrbr,kassel,broadhurst96,bbbl1,bbbl2}
and its generalizations \cite{hoffman99}.
Of those the `depth-length'
shuffle identities (which are also called
`stuffle identities' or `$*$ products' \cite{hoffman97})
 are obviously related, and in fact
equivalent to our `permutation'
identities, as we have convinced ourselves.
Similarly the `weight-length' shuffle identities 
are clearly related to our various 
`partial integration' identities. 
In this case the question of equivalence is more
difficult and will require further investigation.

The reader will have
noted that we did not make use at all
of the precise form of the path integral action.
Our considerations required the
presence of all propagators $g^{(k)}$, as well
as of all the vertices ${\mathcal V}^{p,q}$, however
they did not determine the statistical weights with which
they should appear in the Feynman diagrams.
Our choice of the weights for the propagators
is mainly motivated by the fact that it 
leads to a suggestive form for the free path integral determinant. 
Namely, a simple
application of the ``$\ln\det = \tr\ln$'' identity shows
that, formally,
\footnote{This calculation may be seen as a ``chiral''
generalization of the calculation of the
Scalar QED Euler-Heisenberg Lagrangian 
performed in \cite{ss1}.}
\bear
\ln Z(0,\lambda) &=& {\mathrm const.} +
\sum_{n=1}^{\infty}\lambda^n {\zeta (n)\over n}.
\label{lndetfree}
\ear\no
Comparing this expression with the well-known
formula for the logarithm of the $\Gamma$ function
\bear
\Gamma (1+x) &=&
\exp \biggl[
-\gamma x + \sum_{n=2}^{\infty}{(-1)^n\over n}
\zeta(n) x^n \biggr],
\label{repGamma}
\ear\no
we see that we can identify the free partition function
with the $\Gamma$ function under the assumption
that the ill-defined $\zeta (1)$ appearing in
(\ref{lndetfree}) is renormalized to Euler's constant
$\gamma$,%
\footnote{Considering the identities
$\gamma = \lim_{n\to\infty}\left(\sum_{k=1}^n\frac{1}{k}-\ln n\right)$
and $\sum_{k=1}^{\infty}\frac{1}{k^{\epsilon}} =
{1\over\epsilon} +\gamma + {\rm O}(\epsilon)$
this assumption seems quite natural.}
\bear
Z_{\mathrm renorm}(0,\lambda) &=& {\mathrm const.} \times
\Gamma (1-\lambda).
\label{freepartition}
\ear\no
Similarly the total propagator becomes relatively simply.
Using the integral representation eq.\ (\ref{intrepLik})
of the polylogarithm
it is easily shown that
\bear
p_{12} &\equiv&
\sum_{k=0}^{\infty}
\lambda^k {(2\pi\i)}^kg_{12}^{(k)}
=
g^{(0)}_{12}
+ {z_{12}\lambda\over 1-\lambda}
\phantom{,}_2F_1 (1,1-\lambda;2-\lambda;z_{12}).
\non\\
\label{totprop}
\ear\no

For the interaction term there seems to be no such preferred
choice. A question of obvious interest (but equally obvious
difficulty) is whether non-trivial interaction potentials
$V(g,\bar g)$ exist such that the $\zeta$-model would be
exactly solvable. 

Recently the following generalization of Euler-Zagier
sums (\ref{defmultizeta}) has attracted some attention,
\bear
\zeta(k_1,\ldots,k_m;\sigma_1,\ldots,\sigma_m) &=&
\sum_{n_1>n_2>\cdots >n_m>0} 
{\sigma_1^{n_1}\cdots \sigma_m^{n_m}
\over
n_1^{k_1}\cdots n_m^{k_m}
},
\label{defmultizetaalt}
\ear\no
where $\sigma_j = \pm 1$. Those ``alternating'' Euler-Zagier sums
arise naturally in the calculation of ultraviolet divergences
in renormalizable quantum field theories, and in the application
of knot theory to the classification of those divergences
\cite{brokre96}. At the same time 
the inclusion of alternating Euler-Zagier sums seems to
simplify the problem of reducing the set of all such
sums to a basic set via multiple $\zeta$ identities
\cite{broadhurst96}. (For a tabulation of alternating series see
\cite{bjop}.)
   
More generally, arbitrary $N$-th roots of unity in place of the
$\sigma_j$ have been considered in connection with
the study of mixed Tate motives over Spec ${\mathbb Z}$ 
\cite{goncharov95,racinet}. 
Those phase factors can be easily accomodated in the 
$\zeta$-model. To generate a phase factor $\sigma=\e^{2\pi\i s}$ 
the propagator (\ref{defgk}) has to be simply replaced by
\bear
g^{(k)}_{\sigma}(u_{12})
\equiv
\sum_{n=1}^{\infty}
{\e^{2\pi\i(n+s)u_{12}}
\over
{(2\pi\i n)}^k}.
\label{defgktwist}
\ear\no
On the path integral level this can be achieved
by changing from periodic to twisted boundary
conditions, $x(1) = \sigma x(0)$, and replacing 
$\partial$ by $\partial -2\pi\i s$.

\vskip20pt

\no
{\bf Acknowledgements:}
We would like to thank C. Zeta-Jones for inspiration.
Helpful discussions with
D. Broadhurst and D. Kreimer 
are gratefully acknowledged.
We also thank M. E. Hoffman for detailed comments
on the electronic preprint version of this paper, 
and L. Dixon for pointing out a typo.
One of  us (U. M.) would like to thank the LAPTH for the kind
hospitality during a stay when part of this work was done.

\vfill\eject

\end{document}